\documentclass{article}

\usepackage{amsthm,amsfonts,amssymb,amsmath}
\usepackage{enumerate}
\usepackage[colorlinks=true]{hyperref}

\newtheorem{thm}{Theorem}[section]
\newtheorem{cor}[thm]{Corollary}
\newtheorem{lem}[thm]{Lemma}

\newtheorem{pro}[thm]{Proposition}

\theoremstyle{remark}

\theoremstyle{definition}

\newcommand{\CC}{\mathbb{C}}
\newcommand{\RR}{\mathbb{R}}
\newcommand{\ZZ}{\mathbb{Z}}
\newcommand{\QQ}{\mathbb{Q}}

\newcommand{\DD}{\mathbb{D}}

\newcommand{\CA}{{\mathcal {A}}}

\newcommand{\CE}{{\mathcal {E}}}
\newcommand{\CF}{{\mathcal {F}}}
\newcommand{\CG}{{\mathcal {G}}}

\newcommand{\CL}{{\mathcal {L}}}

\newcommand{\CO}{{\mathcal {O}}}

\newcommand{\CX}{{\mathcal {X}}}

\newcommand{\CZ}{{\mathcal {Z}}}

\newcommand{\an}{{\mathrm{an}}}

\newcommand{\st}{{\mathrm{st}}}

\newcommand{\End}{{\mathrm{End}}}

\newcommand{\Hom}{{\mathrm{Hom}}}

\renewcommand{\Im}{{\mathrm{Im}}}

\newcommand{\Lie}{{\mathrm{Lie}}}

\newcommand{\trdeg}{{\mathrm{trdeg}}}

\newcommand{\Alb}{{\mathrm{Alb}}}

\newcommand{\Sp}{{\mathrm{Sp}}}

\newcommand{\zar}{{\mathrm{zar}}}

\DeclareMathOperator{\Spec}{Spec}

\renewcommand{\d}{\textnormal{d}}

\newcommand{\wt}{\widetilde}

\newcommand{\ds}{\displaystyle}

\newcommand{\ol }{\overline}

\newcommand{\lra}{\longrightarrow}

\newcommand{\sA}{{\mathcal A}}

\newcommand{\C}{{\mathbb C}}
\newcommand{\D}{{\mathbb D}}

\newcommand{\R}{{\mathbb R}}

\begin{document}
%------------------------------------------------------

%\title{The geometric Bombieri--Lang conjecture for varieties of maximal Albanese dimension}

\title{The Geometric Bombieri--Lang Conjecture for Ramified Covers of Abelian Varieties}

\author{Junyi Xie, Xinyi Yuan}

\maketitle
\tableofcontents

\section{Introduction}

The far-reaching Bombieri--Lang conjecture is a high-dimensional generalization of the celebrated Mordell conjecture proved by Faltings \cite{Faltings1983}. 
Beyond the Mordell conjecture, the only known case of the Bombieri--Lang conjecture is 
the case of subvarieties of abelian varieties by Faltings \cite{Faltings1991, Faltings1994}.

The geometric Bombieri--Lang conjecture is an analogue of the Bombieri--Lang conjecture over function fields. We refer to \cite[\S1.1, \S2.4]{XY2023a} for some precise versions of the geometric Bombieri--Lang conjecture in characteristic 0, and we will restrict to characteristic 0 throughout this paper.
The conjecture (or its suitable versions) is proved 
for curves by \cite{Manin63,Grauert65}, for
subvarieties of abelian varieties by  \cite{Raynaud83, Buium92, Hrushovski96},    
for smooth projective varieties with ample cotangent bundles by
\cite{Nog82,MD84}, and for constant hyperbolic varieties by \cite[Cor. 4.2]{Nog92}.

In \cite{XY2023a}, we have proposed an approach to the geometric Bombieri--Lang conjecture for hyperbolic varieties in characteristic 0, which applies the classical Brody lemma to construct entire curves from rational points. The approach is conditional on a non-degeneracy conjecture about partial heights. This non-degeneracy conjecture is confirmed in the case that the variety is finite over an abelian variety, and as a consequence, 
\cite[Thm. 1.2]{XY2023a} confirms the geometric Bombieri--Lang conjecture for hyperbolic varieties finite over abelian varieties.

In this paper,  we move from ``hyperbolic'' to ``general type''. 
The following main theorem confirms the geometric Bombieri--Lang conjecture 
(cf. \cite[Conj. 1.1]{XY2023a}) for varieties finite over abelian varieties of trivial trace.

\begin{thm} [ramified covers of abelian varieties] \label{main}
Let $K$ be a finitely generated field over a field $k$ of characteristic 0. 
Let $X$ be a projective variety over $K$ with a finite morphism $f:X\to A$ for an abelian variety $A$ over $K$. 
Assume that the $\ol K/\bar k$-trace of $A_{\ol K}$ is 0. 
Then $(X\setminus \Sp_{\rm alg}(X))(K')$ is finite for any finite extension $K'$ of $K$.
\end{thm}

Note that we do not assume that $X$ is smooth over $K$ or that $f$ is surjective in the theorem. 
In particular, $X$ is allowed to be a closed subvariety of $A$, and in this case, versions of the geometric Bombieri--Lang conjecture were proved by \cite{Raynaud83, Buium92}. 

We refer to \cite{Conrad} for more details on Chow's trace of abelian varieties. 
The above assumption of trivial $\ol K/\bar k$-trace excludes the complication caused by constant varieties (coming from a base change from $k$ to $K$), so the result becomes very clean. 

Recall that the algebraic special set $\Sp_{\rm alg}(X)$ of $X$ is the Zariski closure in $X$ of the union of the images of all non-constant rational maps from abelian varieties over $\overline K$ to $X_{\overline K}$. 

The Green--Griffiths--Lang conjecture holds for such $X$ in the theorem; i.e., $X$ is of general type if and only if $\Sp_{\rm alg}(X)\neq X$.
Moreover, for any complex variety $Y$ finite over a complex abelian variety,  Lang's conjecture that $\Sp_{\rm alg}(Y)=\Sp_{\rm an}(Y)$ is also confirmed, where the analytic special set $\Sp_{\rm an}(Y)$ of $Y$ is the Zariski closure in $Y$ of the union of all entire curves in $Y$. 
These two results are consequences of works of Ueno, Kawamata and Yamanoi. We refer to \cite[Thm. 4.1, Cor. 4.3]{XY2023a} for more details.

The theorem has the following variant, which also confirms a weak version of \cite[Conj. 1.1]{XY2023a}. 

\begin{cor} [maximal Albanese dimension]
Let $K$ be a finitely generated field over a field $k$ of characteristic 0. 
Let $X$ be a normal projective variety of general type over $K$.
Assume that $X$ has a maximal Albanese dimension, and that the $\ol K/\bar k$-trace of its Albanese variety is 0. 
Then there is a proper Zariski closed subset $Z\subsetneq X$ such that
$(X\setminus Z)(K')$ is finite for any finite extension $K'$ of $K$.
\end{cor}

We refer to \cite[Prop. A.6]{Mochizuki12} for the basics of Albanese varieties. 
The normal projective variety $X$ is said to have a \emph{maximal Albanese dimension} if the Albanese morphism $X_{\ol K}\to \Alb(X_{\ol K})$ (via a base point of $X(\ol K)$) is generically finite onto its image. To prove the corollary, we only need to treat that case $k=\overline{k}.$ After a suitable base change, we may assume that 
 $X(K)\neq \emptyset$ to have an Albanese morphism $X\to \Alb(X)$. Take $X_0$ to be the normalization of 
$\Im(X\to \Alb(X))$ in $X$. This gives a birational morphism $\psi:X\to X_0$. 
Apply Theorem \ref{main} to $X_0$. 
Take $Z$ to be the union of $\psi^{-1}(\Sp_{\rm alg}(X_0))$ with the minimal Zariski closed subset $W$ of $X$ such that $X\setminus W\to X_0\setminus \psi(W)$ is an isomorphism.

In the end, we describe our proof of Theorem \ref{main}.
It is still based on the idea of constructing entire curves in \cite{XY2023a}, but the situation is much more delicate. 
In fact, we first use the Lefschetz principle to reduce the problem to the essential case $k=\CC$ and $K=\CC(B)$ where $B$ is a smooth projective curve over $\C$. 
Take an integral model $\CX\to B$ of $X$ over $B$. 
For the sake of contradiction, assume that there is an infinite sequence $\{x_n\}_{n\geq1}$ in $(X\setminus \Sp_{\rm alg}(X))(K)$. 
As in the proof of \cite[Thm. 1.2]{XY2023a}, we can apply Brody's lemma to obtain an entire curve in a smooth fiber $\CX_b$ for some closed point $b\in B$, but this is not sufficient for our purpose, unless the entire curve is not contained in $\Sp_{\rm alg}(\CX_b)$. 
We need to control the location of the entire curve, and this is usually very hard in the abstract situation of the Brody lemma. 
Our solution of the problem is to use more explicit methods to construct entire curves on fibers of integral models of $A$ over $B$, and ``lift'' the entire curves to fibers of $\CX$ over $B$ by the Brody lemma. 

Let us first introduce a general situation for the abelian variety $A$ over $K$.
Let $Z$ be a Zariski closed subset of $A$. 
Let $\{s_n\}_{n\geq1}$ be an infinite sequence in $(A\setminus Z)(K)$. 
Take an integral model $\CA\to B$ of $A$ over $B$, and denote by $\CZ$ the Zariski closure of $Z$ in $\CA$. 
The goal is to take a limit of a reparametrization of $\{s_n:B\to \CA\}$ to construct an entire curve $\phi:\CC\to \CA_b$ in a smooth fiber $\CA_b$ for some closed point $b\in B$, and ensure that $\phi(\CC)$ is not contained in $\CZ_b$. 
This will be done by explicit construction of entire curves. 
In fact, by careful and precise choices of reparametrization, our limit entire curves are linear entire curves in $\CA_b$ in the sense that the inverse images of the entire curves in the universal covering $\Lie(\CA_b)$ of $\CA_b$ are translations of complex lines in $\Lie(\CA_b)$. 
Note that this directly implies the theorem in the case that 
$f:X\to A$ is a closed immersion.

Now we go back to the original situation.  
By the finite morphism $f:X\to A$, the sequence $\{x_n\}_{n\geq1}$ in $(X\setminus \Sp_{\rm alg}(X))(K)$ gives a sequence $\{s_n\}_{n\geq1}$ in $A(K)$. 
Set $Z=f(\Sp_{\rm alg}(X))$. 
Assume that $f^{-1}(Z)= \Sp_{\rm alg}(X)$, so that
$\{s_n\}_{n\geq1}$ actually lies in $(A\setminus Z)(K)$.
This assumption is not automatic, but it is actually a minor one. 
Take the integral models $\CA\to B$ and the Zariski closure $\CZ$ as above. 
Assume that 
$f$ extends to a morphism $\tilde f: \CX_U\to \CA_U$ for a non-empty open subvariety $U$ of $B$. 
From the above, we have constructed an entire curve $\phi:\CC\to \CA_b$ which is not contained in $\CZ_b$.
Note that $\phi:\CC\to \CA_b$ is obtained as a limit of a reparametrization $\{s_n':\DD\to \CA\}$ of $\{s_n:U\to \CA\}$, 
and thus we can lift the reparametrization $\{s_n':\DD\to \CA\}$ to a reparametrization $\{x_n':\DD\to \CX\}$ of $\{x_n:U\to \CX\}$ by the same style, 
and apply the Brody lemma to $\{x_n':\DD\to \CX\}$ to obtain an entire curve 
$\wt\phi:\CC\to \CX_b$. 
Because of the limit process, the entire curve 
$\wt\phi:\CC\to \CX_b$ is not necessarily a lifting of the entire curve $\phi:\CC\to \CA_b$, but by some explicit consideration, we can prove that 
$\overline{\tilde f(\wt\phi(\CC))}=\overline{\phi(\CC)}$ under the Euclidean topology. 
Then $\tilde f(\wt\phi(\CC)) \not\subseteq \CZ_b$. 
By some choice of $b$, we have $\Sp_{\rm alg}(\CX_b)=\tilde f^{-1}(\CZ_b)$. 
Then we conclude that $\wt\phi(\CC)\not\subseteq\Sp_{\rm alg}(\CX_b)$.
This finishes the idea of the proof.

\medskip

\noindent\textbf{Acknowledgment.}
The authors are grateful to Guoquan Gao, Ariyan Javanpeykar, Junjiro Noguchi, Gang Tian, Songyan Xie,  Chenyang Xu, and Shou-Wu Zhang for many important communications on the first version of this work.

The authors would like to thank the support of the China-Russia Mathematics Center during the preparation of this paper. 
The first author is supported by a grant from the National Science Foundation of China (grant NO.12271007).
The second author is supported by a grant from the National Science Foundation of China (grant NO. 12250004) and the Xplorer Prize from
the New Cornerstone Science Foundation.

\section{Transfer maps and linear entire curves}

In this section, we introduce our notations and prove some preliminary results as preparation for the proof of the main theorem in the next section.
In \S\ref{sec Betti}, we introduce the transfer map between tangent spaces in the setting of the Betti map, and prove a non-degeneracy result on the transfer map. 
In \S\ref{sec linear}, we introduce the notion of linear entire curves, and 
 prove a result on the Zariski closures of linear leaves.

\subsection{Notations and terminology}

In this paper, we resume all the notations and terminology in \cite{XY2023a}. For convenience, we repeat them here. 

For any abelian group $M$ and any ring $R$ containing $\ZZ$, denote
$M_R=M\otimes_\ZZ R$. This apply particularly to $R=\QQ,\RR,\CC$.

By a \emph{variety}, we mean an integral scheme, separated of finite type over the base field. A \emph{curve} is a 1-dimensional variety. 

By a \emph{function field of one variable} over a field $k$, we mean a finitely generated field $K$ over $k$ of transcendence degree 1 such that $k$ is algebraically closed in $K$. 
We usual denote by $B$ a smooth quasi-projective curve over $k$ with function field $K$.
For a projective variety $X$ over $K$, an \emph{integral model of $X$ over $B$}
is a quasi-projective variety $\CX$ over $k$ together with a projective and flat morphism $\CX\to B$ whose generic fiber is isomorphic to $X$.

Let $K$ be a finitely generated field over a field $k$, and assume that $k$ is algebraically closed in $K$.
Let $A$ be an abelian variety over $K$. 
Denote by $A^{(K/k)}$ Chow's $K/k$-trace of $A$ over $k$.
Denote
$$V(A,K):=(A(K)/A^{(K/k)}(k))\otimes_{\ZZ}\RR,$$ 
which is a finite dimensional $\R$-vector space.
In our setting, $k$ is usually $\CC$.

By a \emph{point} of a variety over $\CC$, we mean a closed point. 
By the \emph{generic point} of an integral variety, we mean the generic point of the scheme.

All complex analytic varieties are assumed to be reduced and irreducible.
For a complex analytic variety $X$ with a point $x\in X$, denote by $T_xX$ the complex analytic tangent space of $X$ at $x$ defined by holomorphic derivations. 
For a holomorphic map $f:X\to Y$ with $f(x)=y$, denote by $\d f: T_xX\to T_yY$ the induced map between the tangent spaces.

For a complex abelian variety $A$, the complex Lie algebra $\Lie(A)$ is defined to be the group of translation-invariant holomorphic derivations on 
$\CO_A$. By restriction, we have canonical isomorphisms $\Lie(A)\simeq T_xA$ for any $x\in A$.

Let $X$ be a complex analytic space, and $S$ be a subset of $X$. 
Denote by $\ol S$ the closure of $S$ in $X$ under the Euclidean topology. 
Assume that the smooth locus $X^{\mathrm{sm}}$ is covered by countably many open balls $\{U_\alpha\}_{\alpha\in I}$. 
We say that $S$ is \emph{measurable} if $S\cap D_\alpha$ is measurable under the Lebesgue measure of the ball $D_\alpha$ for all 
$\alpha$.  We say that $S$ has measure zero if $S\cap D_\alpha$ has measurable zero for all $\alpha$.

For any positive real number $r$, denote by
$$
\DD_r:=\{z\in \CC: |z|< r\},\qquad
\overline\DD_r:=\{z\in \CC: |z|\leq r\}
$$
the discs of radius $r$. 
Write $\DD=\DD_1$ and $\overline\DD=\overline\DD_1$.
Denote by $\ds v_{\rm st}=\frac{\d}{\d z}$ the tangent vector of $\DD_r$ at $0$ under the standard coordinate $z$.

Let $(Y, d)$ be a metric space. Let $\{r_n\}_{n\geq 1}$ be a sequence of positive real numbers convergent to infinity. Let $\phi_n:\D_{r_n}\to Y$ be a sequence of continuous maps.
We say that $\{\phi_n\}_n$ converges to a map $\phi: \C\to Y$ if it converges on every compact subset $\Omega$ of $\C$. Since $\Omega\subseteq \D_{r_n}$ for $n$ sufficiently large, the above definition makes sense. If such $\phi$ exists, it is unique and continuous. Moreover, if $Y$ is further a complex analytic variety and $\phi_n$ are holomorphic, then $\phi$ is holomorphic. We say that $\phi_n$ \emph{uniformly} converges  to $\phi: \C\to Y$ if for every $\epsilon>0$, there is $N\geq 1$ such that for every $n\geq N$ and $z\in \D_{r_n}$, $d(\phi(z),\phi_n(z))\leq \epsilon.$ If $\phi_n$ uniformly converges to $\phi$, then it converges to $\phi.$ 
%One has the following Cauchy convergence criterion:
%$\phi_n$ uniformly converges to some $\phi$ if and only if for every $\epsilon>0$ there is $N\geq 1$ such that for every $m,n\geq N$ and $z\in \D_{\min\{r_n.r_m\}}$ 
%$d(\phi_m(z),\phi_n(z))\leq \epsilon.$

\subsection{Betti maps and transfer maps} \label{sec Betti}

The Betti maps for complex abelian schemes were first introduced by Mok \cite[p. 374]{Mok1989}. They were reviewed in \cite[\S2]{Cantat2021} and \cite[\S3.1]{XY2023a}. 
We will recall some related notations in \cite[\S3.1]{XY2023a}, and introduce a new term called the transfer map.

\subsubsection*{Betti maps}

Let $B$ be a Riemann surface.
Let $\pi:\CA\to B$ be a holomorphic family of abelian varieties over $B$; that is, $\CA$ is a complex manifold, $\pi$ is a smooth holomorphic map endowed with a holomorphic section $e:B\to\CA$, and every fiber of $\pi$ is an abelian variety with the identity point induced by $e$. 
Denote by $\CA(B)_h$ the group of holomorphic sections $s:B\to \CA$. 

If $B, \CA, \pi, e$ are algebraic, we will denote by 
$\CA(B)$ the group of algebraic sections $s:B\to \CA$. 
Then we have a natural injection $\CA(B)\to \CA(B)_h$.
However, in this subsection, we take the more general analytic setting.

Let $U\subset B$ be a connected and \emph{simply connected} open neighbourhood of $b$ in $B$. 
Let $b\in B$ be a point.
Denote by $\pi_U:\CA_U\to U$ the base change of $\pi$, and by $\CA_b$ the fiber of $\CA$ above $b$.
The \emph{Betti map} is a canonical real analytic map 
$$
\beta=\beta_{b,U}: \CA_U \lra \CA_b,
$$
which satisfies the following properties:
\begin{enumerate}[(1)]
\item The composition 
$\CA_{b}\hookrightarrow\CA_U\stackrel{\beta}{\to} \CA_b$ is the identity map. 
\item For any point $b'\in U$, the composition 
$\CA_{b'}\hookrightarrow\CA_U\stackrel{\beta}{\to} \CA_b$ is an isomorphism of real Lie groups.
\item The induced map 
$$\tilde\beta=(\beta,\pi): \CA_U \lra \CA_b\times U$$
 is a real analytic diffeomorphism of manifolds.
\item For any $x\in \CA_U$, the fiber $\beta^{-1}(\beta(x))$ is a complex analytic subset of $\CA_U$, which is biholomorphic to $U$. 
\end{enumerate}

The fiber $\CF_{x,U}=\beta_{b,U}^{-1}(\beta_{b,U}(x))$ is called \emph{the local Betti leaf at $x$ over $U$}. 
It is independent of the choice of $b\in U$ (for fixed $x,U$), and its 
germ $\CF_x$ at $x$ is independent of the choice of $U$. 
A \emph{Betti leaf} of $\CA$ is a connected subset $\CF_0$ of $\CA$ such that for any $x\in \CF_0$, the germ of $\CF_0$ at $x$ is equal to $\CF_x$. 
Note that any connected component of a torsion multi-section of $\CA\to B$ is a Betti leaf. 
The set of all Betti leaves forms a \emph{Betti foliation}.

\subsubsection*{Transfer maps between tangent spaces}

For any $x\in \CA$, there is a direct sum decomposition
$$
T_{x}\CA = T_x \CA_{\pi(x)}\oplus T_x \CF_x 
$$
of complex analytic tangent spaces at $x$. 
Here $\CF_x$ is the germ of a local Betti leaf through $x$.
We call the induced homomorphism 
$p_1:T_{x}\CA \to T_x \CA_{\pi(x)}$
the \emph{Betti projection}.

The map $p_1:T_{x}\CA \to T_x \CA_{\pi(x)}$ is equal to the map $\d \beta:T_{x}\CA \to T_x \CA_{\pi(x)}$ induced by a Betti map $\beta:\CA_U\to \CA_{\pi(x)}$. Then $\d \beta$ is a $\CC$-linear map, which is a priori only $\RR$-linear.

Let $s:B\to \CA$ be a section of $\pi:\CA\to B$. 
For any point $b\in B$, define \emph{the transfer map}
$$
\delta(s, \cdot)=\delta_b(s, \cdot): T_b B \lra \Lie(\CA_b)  
$$
to be the composition 
$$
T_b B \stackrel{\d s}{\lra} T_{s(b)}(\CA)  \stackrel{p_1}{\lra} T_{s(b)}(\CA_b)
\lra \Lie(\CA_b). 
$$
Here the last map is the canonical isomorphism via translation-invariant vector fields.
We have the following basic properties.

\begin{lem}\label{transfer tangent}
\begin{enumerate}[(1)]
\item The transfer map $\delta(s, \cdot):T_b B \to \Lie(\CA_b)$ is additive in $s\in \CA(B)_h$.
\item Let $A_0$ be a complex abelian variety with a homomorphism $i:A_0\times B\to \CA$ over $B$. 
Then $\delta(i(s), \cdot):T_b B \to \Lie(\CA_b)$ is 0 for any $s\in A_0$.
Here $i(s)$ denotes the image of the natural map $A_0\to \CA(B)_h$ induced by $i$. 
\end{enumerate}
\end{lem}

\begin{proof}

For (1), it suffices to check it in the setting of differentiable manifolds and real Lie groups. Take an open disc $U$ of $B$ containing $x$.
Then $\delta(s, \cdot)$ is the map between the tangent spaces induced by 
$\beta\circ s$, i.e. the composition
$$
U \stackrel{s}{\lra} \CA_U \stackrel{\beta}{\lra} \CA_b.
$$
Note that $\beta$ is additive. 
Then the additivity follows from the classical statement that for two smooth maps 
$f_1, f_2: M\to G$ from a manifold $M$ to a Lie group $G$, the induced maps $\d f_i: T_mM\to \Lie(G)$
at a point $m\in M$ are additive in that $\d(f_1+f_2)=\d f_1+\d f_2$.

For (2), it suffices to note that the image of the composition $\beta\circ i(s):U \to \CA_b$ is a point of $\CA_b$. 
This follows from the fact that the composition factors through a similar map for the trivial family $A_0\times B\to B$.
\end{proof}

In the setting of the lemma, by continuity, the transfer map induces a canonical map 
$$
\delta(s, \cdot): T_b B \lra \Lie(\CA_b)  
$$
for any $s$ in $(\CA(B)_h/A_0)_\RR=(\CA(B)_h/A_0)\otimes_\ZZ \RR$. 
This can also be written as a bilinear map 
$$
\delta(\cdot, \cdot): (\CA(B)_h/A_0)_ \RR\times T_b B \lra \Lie(\CA_b).
$$

\subsubsection*{Non-degeneracy of the transfer maps} 

Now we restrict to the algebraic situation. 
The treatment of the non-degeneracy result in \cite[\S3.2]{XY2023a} has the following important consequence, which will be crucially used in the proof of Theorem \ref{main}. 

\begin{thm}[non-degeneracy] \label{non-deg2}
Let $B$ be a complex smooth quasi-projective curve, and let $\pi:\CA\to B$ be an abelian scheme.
Let $K=\CC(B)$ be the function field, and let $A=\CA_K$ be the generic fiber of $\CA$.  
Let $s\in (A(K)/A^{(K/\CC)}(\CC))\otimes_\ZZ\RR$ be a nonzero element. 
Let $\Sigma(s)$ be the set of $b\in B$ such that the transfer map
$$\delta(s,\cdot):T_bB\lra \Lie(\CA_b)$$
is  zero.
Then $\Sigma(s)$ has measure 0 in $B$ with respect to any K\"ahler form on $B$. 
\end{thm}

\begin{proof}

Let us first introduce some extra notations as in \cite[\S3.2]{XY2023a}.
Take a symmetric and ample line bundle $L$ on $A$, and extend it to a line bundle $\CL$ on $\CA$. 
Let $\omega=\omega(\CL)$ be the Betti form  on $\CA$ associated to $\CL$ as introduced in \cite[\S3.1]{XY2023a}. 
By Gauthier--Vigny \cite[Thm. B]{GV2019} (cf. \cite[Thm. 3.1]{XY2023a}), for any section $s\in \CA(B)\simeq A(K)$, the canonical height 
$$
\hat h_{L}(s)=\int_{B} s^*\omega.
$$

Now we are ready to prove the current theorem. 
To illustrate the idea, we first consider the case $s \in A(K)$ with $\hat h_L(s)\neq 0$.
Denote by $\Sigma(s)'$ the set of $b\in B$ such that the fiber $(s^*\omega)(b)$ in 
$(T_bB)^\vee \otimes_\CC \overline{(T_bB)}^\vee$ is 0. 
The result is the consequence of the following two properties:
\begin{enumerate}[(a)]
\item $\Sigma(s)\subset\Sigma(s)'$;
\item $\Sigma(s)'$ has measure zero in $B$.  
\end{enumerate}

We first prove (b).
By $\displaystyle \hat h_L(s)>0$, we see that 
$s^*\omega$ is not identically zero on $B$.
As in the proof of \cite[Thm. 3.5]{XY2023a}, 
 $\Sigma(s)'\subsetneq B$ is the zero locus of a real analytic function on $B$, so $\Sigma(s)'$ has measure zero in $B$. 

Now we prove (a).
Recall the Betti map $\beta:\CA \stackrel{}{\to} \CA_b$
and the Betti form $\omega=\beta^*\omega_b$.
Here $\omega_b$ is a translation-invariant positive $(1,1)$-form on $\CA_b$ representing $c_1(\CL_b)$. 
By translation, we can write $\omega_b$ as a finite sum
$$\omega_b=\sum_{j=1}^d \gamma_j\wedge \bar\gamma_j',$$
where $\gamma_j$ and $\gamma_j'$ are translation-invariant holomorphic 1-forms on $\CA_b$. 
We further have
$$
s^*\omega= (\beta\circ s)^*\omega_b= \sum_{j=1}^d (\beta\circ s)^*\gamma_j\wedge (\beta\circ s)^*\bar\gamma_j'.
$$

Recall that
$$\delta(s,b):T_bB\lra \Lie(\CA_b)$$
is exactly the push-forward via the composition $\beta\circ s: B \to \CA_b$. 
Then this push-forward is zero for any $b\in \Sigma(s)$.
Take a dual, we see that $((\beta\circ s)^*\gamma_j)(b)=0$ for any $j$. 
It follows that $(s^*\omega)(b)=0$, and thus $b\in \Sigma(s)'$.

Now we extend the above proof to any nonzero 
$s\in (A(K)/A^{(K/\CC)}(\CC))\otimes_\ZZ\RR$. 
Note that $\displaystyle \hat h_L(s)>0$ by the Lang--N\'eron theorem (cf. \cite[Thm. 3.4]{XY2023a}). 
We only need to explain the term ``$s^*\omega$'' and the term ``$(\beta\circ s)^*\gamma_j$'' satisfying the above properties. 
The term $s^*\omega$ is given by the $(1,1)$-form $\omega(s,s)$ constructed in 
\cite[Lem. 3.3]{XY2023a}. 
For the term $(\beta\circ s)^*\gamma_j$, it suffices to define $(\beta\circ s)^*: \Lie(\CA_b)^\vee \to (T_bB)^\vee$
as the dual of $\delta(s,b):T_bB\to \Lie(\CA_b)$. 
The relation between $s^*\omega$ and $(\beta\circ s)^*\gamma_j$ is extended from 
$A(K)$ to $(A(K)/A^{(K/\CC)}(\CC))\otimes_\ZZ\RR$ by linearity and a limit process. 
This finishes the proof. 
\end{proof}

\subsection{Linear entire curves} \label{sec linear}

The goal of this section is to introduce linear entire curves and study some basic properties.

\subsubsection*{Linear entire curves}

Let $A$ be a complex abelian variety. 
Let $x\in A$ be a point, and $v\in \Lie(A)$ be a nonzero vector in the Lie algebra.
By definition, $v$ defines a translation-invariant vector field on $A$, and thus it is a foliation $\CF(v)$ on $A$. 
The leaf $\CF(v)_x=\CF(x,v)$ of this foliation through $x$ is called 
the \emph{linear leaf through $x$ in the direction $v$}. 
By translation, we have
$$
\CF(x,v)=\CF(0,v)+x.
$$

Alternatively, write $A=\Lie(A)/H_1(A,\ZZ)$. 
Denote by $\tilde x$ a lifting of $x$ in $\Lie(A)$.
In the complex vector space $\Lie(A)$,  there is a linear subset given by 
$$
\wt\phi_{(\tilde x,v)}:\CC \lra \Lie(A), \quad z \longmapsto \tilde x+z v.
$$
Composing with $\Lie(A)\to A$, we obtain an entire curve
$$
\phi_{(x,v)}:\CC \lra A.
$$
The entire curve is independent of the choice of $\tilde x$, and its image is exactly the linear leaf $\CF(x,v)$. 
We call $\phi_{(x,v)}:\CC\to A$ \emph{the linear entire curve through $x$
in the direction $v$}.

The linear leaf $\CF(x,v)$ is generally not algebraic in $A$.
Its Zariski closure is equal to $A(v)+x$ for an abelian subvariety $A(v)$ of $A$.
In fact, by translation, it suffices to assume $x=0$. 
Then $\phi_{(0,v)}$ is a homomorphism, so its image $\CF(0,v)$ is actually a complex Lie group, and thus its Zariski closure $A(v)$ inherits a group structure.
Moreover, $A(v)$ is irreducible as an analytic variety, since the image of $\CC$ 
(under $\phi_{(0,v)}$) is irreducible. 
Then $A(v)$ is an abelian subvariety of $A$.

\subsubsection*{Specialization of endomorphism ring}

The following result will be used later. It is well-known, but we include a quick proof in terms of the concept of transcendental closed points in \cite[\S4.1]{XY2023a}.

\begin{lem} \label{transcendental2}
Let $B$ be a smooth curve over $\CC$, and $\pi:\CA\to B$ be an abelian scheme over $B$.
Denote $K=\CC(B)$, so $A=\CA_K$ is an abelian variety over $K$. 
Assume that the natural map 
$\End_K(A)\to \End_{\ol K}(A_{\ol K})$
is an isomorphism. 
Then there are only countably many closed points $b\in B(\CC)$
such that the canonical map 
$\End_K(A)\to \End_\CC(\CA_b)$
induced by specialization is not an isomorphism. 
\end{lem}
\proof
There is a finitely generated subfield $k_0$ of $\CC$ such that 
$(B, \CA, \pi)$ is the base change of a triple $(B_0, \CA_0, \pi_0)$ from $k_0$ to $\CC$ such that the natural map 
$\End_{K_0}(A_0)\to \End_{K}(A_{K})$
is an isomorphism. Here we denote $K_0=k_0(B_0)$ and $A_0=\CA_{0,K_0}$ as before. 
In fact, it suffices to note that $\End_{B}(\CA)=\End_{K}(A_{K})$ is a finitely generated abelian group, so we only need to descend finitely many morphisms of $\End_{B}(\CA)$ to $\End_{B_0}(\CA_0)$ to have 
$\End_{B_0}(\CA_0)\simeq \End_{B}(\CA)$.

We claim that for any {transcendental closed point $b\in B(\CC)$ with respect to $B_0/k_0$}, the canonical map 
$\End_K(X)\to \End_\CC(\CX_b)$ is an isomorphism. 
If so, the result follows from \cite[Lemma 4.4]{XY2023a}.

For the claim, it suffices to check that the map 
$\End_{K_0}(A_0)\to \End_\CC(A_{0,\CC})$ induced by any injection 
$\iota:K_0\to \CC$ is an isomorphism. 
By assumption, the composition
$\End_{K_0}(A_0)\to \End_K(A)\to \End_{\ol K}(A_{\ol K})$
is an isomorphism, 
so $\End_{K_0}(A_0)\to \End_{\ol K_0}(A_{0,\ol K_0})$
is an isomorphism. 
Then it suffices to check that the map
$\End_{\ol K_0}(A_{0,\ol K_0})\to \End_\CC(A_{0,\CC})$ induced by 
any injection $\ol K_0\to \CC$ 
is an isomorphism.
This follows a well-known theorem due to Chow (cf. \cite[Thm. 3.19]{Conrad}). 
\endproof

\subsubsection*{Zariski closures of linear leaves}

The main result of this subsection is to give a clear description of Zariski closures of linear leaves in fibers of abelian schemes.

Let $B$ be a smooth quasi-projective curve over $\CC$, and let $\pi:\CA\to B$ be an abelian scheme.
Let $K=\CC(B)$ be the function field of $B$, and let $A$ be the generic fiber of $\pi:\CA\to B$. 
By the Lang--N\'eron theorem (cf. \cite[Thm. 3.4]{XY2023a}), 
$$V(A,K)=(A(K)/A^{(K/\CC)}(\CC))_\RR$$ 
is a finite dimensional $\R$-vector space.
Note that an isogeny $A\to A'$ induces an isomorphism $V(A,K)\to V(A',K)$. 
For an abelian subvariety $H$ of $A$,  we may check that $$V(H,K)=(H(K)/H^{(K/\CC)}(\CC))_\RR=(H(K)/(H(K)\cap A^{(K/\CC))}(\CC))_\RR.$$
The natural map $V(H,K)\subseteq V(A,K)$ is injective,
as a consequence of the fact that $A$ is isogenous to $H\times H'$ for an abelian variety $H'$ over $K$.

%Let $H_1,H_2$ be two abelian subvarieties of $A$ with $H_1+H_2=A$ and $\dim H_1+\dim H_2=\dim A$, the morphism $V(H_1,K)\oplus V(H_2, K)\to V(A,K)$ sending $(s_1,s_2)$ to $s_1+s_2$ is an isomorphism. 
% Indeed, the projection $\pi_1:V(A, K)\to V(H_1,K)$ is defined in the following way: For every $x\in A(K)$, $(x+H_2)\cap H_1$ is finite and non-empty. Pick $y\in ((x+H_2)\cap H_1)(\overline{K})$, $((x+H_2)\cap H_1)(\overline{K})=y+(H_1\cap H_2)(\overline{K}).$ There is $m\in \Z_{\geq 1}$ such that $[m](H_2)\cap H_1)=0.$ Then $[m](y)\in H_1(K)$ and does not depend on the choice of $y.$ Then $\pi_1(x)$ is the image of $m^{-1}([m]y)$ in $V(H_1,K).$

\begin{lem}\label{lemintersectionsubgv}Let $H_1,H_2$ be two abelian subvarieties of $A$, and let $H$ be the identity component of $H_1\cap H_2$.
Then $V(H_1, K)\cap V(H_2, K)=V(H, K).$
\end{lem}
\proof 
We may assume that $A=H_1+H_2.$
The up to isogeny, we can decompose 
$$
H_1\sim H\times H_1',\quad
H_2\sim H\times H_2',\quad 
A\sim H\times H_1'\times H_2'.
$$
Then the result is easy.
\endproof

As a consequence, for every $s\in V(A, K)$, there exists a \emph{unique minimal abelian subvariety} $G(s)$ of $A$ with $s\in V(G(s),K).$
A basic property of $G(s)$ is that it is compatible with homomorphisms.

\begin{lem}\label{lempreimageva}
Let $f:A\to A'$ be a homomorphism of abelian varieties over $K$.
Denote by $\bar f: V(A,K)\to V(A',K)$ the map induced by $f$. Then for every 
$s\in V(A,K)$, we have $G(\bar f(s))=f(G(s))$ in $A'$. 
\end{lem}

\proof
Replacing $A'$ by $f(A)$ if necessary, we can assume that $f$ is surjective. 
As $s\in V(G(s),K)$, we have
$\bar f(s)\in V(f(G(s)),K)$, so $G(\bar f(s))\subseteq f(G(s))$ by the minimality of $G(\bar f(s))$.
It remain to prove $f(G(s))\subseteq G(\bar f(s))$.

As we have assumed that $f$ is surjective, $A$ is isogenous to $A'\times A''$ for an abelian subvariety $A''$ of $A$.
Via $V(A,K)= V(A',K)\oplus V(A'',K)$, write $s=s'\oplus s''$. 
Then we have $G(s)\subseteq G(s')\times G(s'')$ by the minimality of $G(s)$.
Apply $f$ to this inclusion, we have $f(G(s))\subseteq G(s')$. 
This finishes the proof.
\endproof

Recall that the transfer map 
$\delta(s,\cdot):T_{b}B\to \Lie(\CA_{b})$ is defined in \S\ref{sec Betti}. 
The following is the main result of this subsection.

\begin{pro} \label{equality of zariski closure}
Let $B$ be a smooth quasi-projective curve over $\CC$, and let $\pi:\CA\to B$ be an abelian scheme.
Let $K=\CC(B)$ be the function field of $B$, and let $A$ be the generic fiber of $\pi:\CA\to B$. 
Assume that the natural map 
$\End_K(A)\to \End_{\ol K}(A_{\ol K})$
is an isomorphism. 
Let $s$ be a nonzero element of $V(A,K)$. 
%Denote by $G$ the smallest abelian subvariety of $A$ such that 
%$s$ lies in $V(G,K)$. 
Denote by $\CG$ the Zariski closure of $G(s)$ in $\CA$. 
Then there is a subset $S(A)$ of $B(\CC)$ of measure 0 such that for any $b\in B(\CC)\setminus S(A)$, the Zariski closure of the linear leaf 
$\CF(0,\delta(s,v_b))$ in $\CA_b$ 
is equal to $\CG_b$. Here $v_b\in T_bB$ is a nonzero vector.
\end{pro}

%\subsubsection*{Zariski closures of linear leaves}
%The following result gives a clear description of Zariski closures of linear leaves in fibers of abelian schemes. Recall that the transfer map 
%$\delta(s,\cdot):T_{b}B\to \Lie(\CA_{b})$ is defined in \S\ref{sec Betti}. 
%
%
%\begin{lem} \label{equality of zariski closure}
%Let $B$ be a smooth quasi-projective curve over $\CC$, and let $\pi:\CA\to B$ be an abelian scheme.
%Let $K=\CC(B)$ be the function field of $B$, and let $A$ be the generic fiber of $\pi:\CA\to B$. 
%Assume that the natural map 
%$\End_K(A)\to \End_{\ol K}(A_{\ol K})$
%is an isomorphism. 
%Let $s$ be a nonzero element of $(A(K)/A^{(K/\CC)}(\CC))_\RR$. 
%Denote by $G$ the smallest abelian subvariety of $A$ such that 
%$s$ lies in $(G(K)/G^{(K/\CC)}(\CC))_\RR$. 
%Denote by $\CG$ the Zariski closure of $G$ in $\CA$. 
%Then there is a subset $S(A)$ of $B(\CC)$ of measure 0 such that for any $b\in B(\CC)\setminus S(A)$, the Zariski closure of the linear leaf 
%$\CF(0,\delta(s,v_b))$ in $\CA_b$ 
%is equal to $\CG_b$. Here $v_b\in T_bB$ is a nonzero vector.
%\end{lem}

\begin{proof}

For any $b\in B(\CC)$, denote by 
$H(A)_b$ the Zariski closure of $\CF(0,\delta(s,v_b))$ in 
$\CA_b$. As $s$ lies in $V(G,K)$,  
$\delta(s,v_b)$ lies in $\Lie(\CG_b)\subseteq \Lie(\CA_b)$. 
It follows that $H(A)_b\subseteq \CG_b$. 
Therefore, we can replace $A$ by $G(s)$ in the lemma. 
As a consequence, we can assume that
$A=G(s)$.
The goal is to prove $H(A)_b=\CA_b$ for $b\notin S$ in this case.

We prove the lemma for by induction on the number of simple isogeny factors of $A$. 
If $A$ is simple, the result is automatic by Lemma \ref{transcendental2} and 
Theorem \ref{non-deg2}. 
In general, denote by $S_0$ the subset of $b\in B(\CC)$ such that the canonical map $\End_K(A)\to \End_\CC(\CA_b)$
is not an isomorphism. 
By Lemma \ref{transcendental2}, $S_0$ has measure 0. 

Let $b\in B(\CC)\setminus S_0$ be an element such that $H(A)_b\neq \CA_b$. 
Then 
$H(A)_b/N\neq \CA_b/N$ for some simple abelian subvariety $N$ of $H(A)_b$. 
As $\End_K(A)\simeq \End_\CC(\CA_b)$, there is a simple abelian subscheme $\CA'$ of $\CA$ such that $N=\CA'_b$. 
Thus $H(A)_b/\CA'_b\neq \CA_b/\CA_b'$ for some simple abelian subvariety $A'$ of $A$, where $\CA'$ denotes the Zariski closure of $A'$ in $\CA$.  
By Lemma \ref{lempreimageva}, $A/A'=G(s')$ where $s'$ is the image of $s$ in $V(A/A',K).$  
Note that 
$H(A/A')_b$ is the image of $H(A)_b$ in $(\CA/\CA')_b$.
It follows that 
$H(A/A')_b\neq (\CA/\CA')_b$, 
and thus $b\in S(A/A')$. 
By induction, the subset $S(A/A')$ of $B(\CC)$ has measure 0. 
As $\End_K(A)$ is a free abelian group of finite rank, 
there are only countably many simple abelian subvarieties $A'$ of $A$. 
Taking the union of $S(A/A')$ over all such $A'$, we still get a subset of $B(\CC)$ of measure 0. 
This finishes the proof.
\end{proof}

\section{Construction of linear entire curves} \label{sec general}

In this section, we prove our main theorem (Theorem \ref{main}).
We will first prove some technical results on linear entire curves as limits of re-parametrizations of sections of abelian schemes in \S\ref{sec limit} and \S\ref{sec location}, and then prove Theorem \ref{main} in \S\ref{sec proof}.

\subsection{Linear entire curves as limits} \label{sec limit}

The major goal of this subsection is to prove the following technical theorem. It gives a concrete construction and a precise description of  linear entire curves obtained from sections of abelian schemes by a Brody-type convergence.

\begin{thm} \label{linear entire curve}
Let $B$ be a smooth quasi-projective curve over $\CC$, and let $\pi:\CA\to B$ be an abelian scheme.
Endow $\CA$  with a K\"ahler metric $\alpha_{\CA}$. 
Let $K=\CC(B)$ be the function field of $B$, and $A$ be the generic fiber of $\pi:\CA\to B$. 
Let $\{s_n\}_{n\geq1}$ be a sequence in $\CA(B)$ and $\{b_n\}_{n\geq1}$ be a sequence in $B$ satisfying the following properties:
\begin{enumerate}[(1)]
\item $b_n$ converges to a point $b\in B$, and 
$s_n(b_n)$ converges to a point $x\in\CA_{b}$.
\item There is a sequence $\ell_n$ of positive real numbers converging to infinity such that the image of $\ell_n^{-1}s_n$ in 
$V(A,K)$
converges to a nonzero element $s_\infty$ in $V(A,K)$. 
\item The transfer map $\delta(s_\infty,\cdot):T_{b}B\to \Lie(\CA_{b})$ is nonzero. 
\end{enumerate}
Define a re-parametrization $\{\phi_n:\DD_{r_n}\to \CA\}_{n\geq1}$ of 
$\{s_n:B\to \CA\}_{n\geq1}$ as follows.
Take an open unit disc $\DD$ in $B$ with center $b$, and assume that $b_n\in \DD$ for all $n\geq 1$. 
Let $z$ be the standard coordinate of $\DD$. Take the tangent vector $v_b\in T_bB$ to be the one represented by $v_\st=\ds\frac{\d}{\d z}\in T_0\DD$. 
Define the re-parametrization map
$$\phi_n: \DD_{r_n}\lra \CA, \quad z\longmapsto s_n(b_n+\ell_n^{-1}z)
$$ 
Here the sum $b_n+\ell_n^{-1}z$ is taken in $\DD$, and 
$\{r_n\}_{n\geq1}$ is a sequence of positive numbers satisfying 
$$r_n/\ell_n<1-|b_n|, \quad 
r_n\lra \infty.$$ 
Then $\{\phi_n:\DD_{r_n}\to \CA\}_{n\geq1}$ satisfies the following properties:
\begin{enumerate}[(a)]
\item The linear entire curve 
$\phi_{(x,\delta(s_\infty,v_b))}$ of $\CA_b$ through $x$ in the direction of $\delta(s_\infty,v_b)$ , viewed as an entire curve of $\CA$ via the embedding $\CA_b\to \CA$, is the unique limit of $\{\phi_n\}_n$. 
\item With the  K\"ahler metric $\alpha_{\CA}$ on
$\CA$,
$$\inf_n\|(\d \phi_n)_{0}(v_{\rm st})\|_{\alpha_{\CA}}>0.$$
\item There exists a choice of the sequence $\{r_n\}_{n\geq1}$ 
(to define $\{\phi_n\}_{n\geq1}$), depending on $\{s_n\}_{n\geq1}$ and  $\{b_n\}_{n\geq1}$, such that 
for every entire curve $\psi:\CC\to \CA$ which is a limit of a re-parametrization of $\{\phi_n:\DD_{r_n}\to \CA\}_{n\geq1}$, the image
$\psi(\CC)$ is contained in 
$\CF(y,\delta(s_\infty,v_b))$ for some $y\in \overline{\CF(x,\delta(s_\infty,v_b))}
\subseteq \CA_b$, and the complement $\CF(y,\delta(s_\infty,v_b))\setminus \psi(\CC)$ has at most one element. 
Here $\overline\Omega$ denotes the closure of a set $\Omega$ under the Euclidean topology.
As a consequence, 
 $$\overline{\psi(\C)}=\overline{\CF(y,\delta(s_\infty,v_b))}=\overline{\CF(x,\delta(s_\infty,v_b))}.$$
\end{enumerate}
\end{thm}

\medskip
For the conditions of the theorem, note that (2) implies that 
$\hat h_L(s_n)$ 
converges to infinity (for a symmetric and ample line bundle $L$).
Conversely, if $\hat h_L(s_n)$ converges to infinity, then (1) and (2) can be satisfied by passing to a subsequence of $\{s_n\}_{n\geq1}$, and (3) can be satisfied for $b\in B$ outside a set of measure 0 by Theorem \ref{non-deg2}. 
This is easy for (1) by the compactness of $\CA_b$. For (2), it comes from the Lang--N\'eron theorem (cf. \cite[Thm. 3.4]{XY2023a}). 
Moreover, we can actually take $\ell_n=\hat h_L(s_n)^{1/2}$ in (2).

The proof of Theorem \ref{linear entire curve} will take up the rest of this subsection. We first introduce some extra notations related to the construction. 

Recall the open disc $\DD$ in $B$ with center $b$. 
%Let $z$ be the standard coordinate of $\DD$. Take the tangent vector $v_b\in T_bB$ to be the one represented by $\ds\frac{\d}{\d z}$. 
%Replacing $\DD$ by a sub-disc if necessary, we assume that $\ol\DD$ is also contained in $B$. 
By abuse of notations, we identify points of $\DD$ with the corresponding one in $B$. In particular, 
we identify $b\in B$ with $0\in \DD$ and thus $\CA_b$ with $\CA_0$. 
Define a map
$$
p_n: \DD_{r_n} \lra \DD, \quad z\longmapsto b_n+\ell_n^{-1}z.
$$
Then the re-parametrization map can be written as
$$\phi_n=s_n|_{\DD}\circ p_n: \DD_{r_n}\lra \CA_\DD, \quad z\longmapsto s_n(b_n+\ell_n^{-1}z)
$$

\subsubsection*{Proof of Theorem \ref{linear entire curve}(b)}
Denote by $v_{b_n}\in T_{b_n}B$ the tangent vector $\ds\frac{\d}{\d z}$ with respect to the coordinate $z$ of $\DD$. 
Take the Betti form $\omega$ on $\CA$ associated to a symmetric and relative ample line bundle $\CL$ on $\CA$.  
By compactness, there is a constant $c>0$ such that 
$\alpha_{\CA} -c\cdot \omega$ is positive definite everywhere on an neighborhood of $\CA_b$.
To prove property (b),  it suffices to prove
$$\inf_n\|(\d \phi_n)_0(v_{\rm st})\|_{\omega}>0.$$
Note that the Betti map $\beta:\CA_\DD\to \CA_b$ does not change the semi-norm of tangent vectors with respect to $\omega$. 
It follows that 
$$ 
\|(\d \phi_n)_0(v_{\rm st})\|_{\omega}
=\ell_n^{-1}\|(\d s_n)_{b_n}(v_{b_n})\|_{\omega}
= \ell_n^{-1}\| \delta_{b_n}(s_n,v_{b_n})\|_{\omega}.
$$
We claim that the last term converges to
$\| \delta(s_\infty,v_{b})\|_{\omega}$, and thus is strictly positive.

Now we prove the claim.
For the case $b_n=b$ for every $n$, we simply have  $\ell_n^{-1} \delta(s_n,v_b)$ converges to 
$\delta(s_\infty,v_b)$ by the $\RR$-linearity of $\delta(s,\cdot)$ in $s$.
In general, the Betti map induces an isomorphism 
$(\d \beta)_{b_n,b}: \Lie(\CA_{b_n})\to \Lie(\CA_{b})$. 
Then it suffices to prove that 
$ (\d \beta)_{b_n,b}(\delta_{b_n}(\ell_n^{-1} s_n,v_{b_n}))$
converges to 
$\delta(s_\infty,v_b)$
in $\Lie(\CA_{b})$. 
By the $\RR$-linearity of $\delta(s,\cdot)$ of $s$ in the finite-dimensional $\RR$-space $(A(K)/A^{(K/\CC)}(\CC))_\RR$, it is reduced to prove the case that $s_n=s_\infty\in \CA(B)$ for every $n$. 
This case is proved following the definitions. 
This proves Theorem \ref{linear entire curve}(b).
%the property (b). 

\subsubsection*{Estimation of errors}

Here we introduce a basic estimate which will be used in the proof of Theorem \ref{linear entire curve}(a)(c).  The idea of the estimate is simply a Taylor expansion, but its writing is rather involved due to various notations of the situation.

Recall that the Betti map gives a real analytic isomorphism
$$\tilde\beta=(\beta,\pi): \CA_\DD \lra \CA_0\times \DD.$$
We fix an euclidean norm $|\cdot|$ on $\Lie(\CA_0)$.
The covering map $\Lie(\CA_0)\to \CA_0$
 induces a distance function $d_{\CA_0}(\cdot,\cdot)$ on $\CA_0$ by 
$$d_{\CA_0}(x,y):=\min\{|x^{\dagger}-y^{\dagger}|\}$$ 
where $x^{\dagger}, y^{\dagger}$ are taken over all liftings of $x$ and $y$ in $\Lie(\CA_0)$ respectively.  
It induces a distance function $d(\cdot,\cdot)$ on $\CA_0\times \DD$ by 
$$d((x_1,t_1),(x_2,t_2)):= d_{\CA_0}(x_1,x_2)+ |t_1-t_2|.$$
Via the real analytic isomorphism  $\tilde\beta$, it induces a distance function on $\CA_\DD$. We still denote the induced distance function by  $d(\cdot,\cdot)$.

Now we have the following estimate.
\begin{lem} \label{error}
There is a positive integer $n_0$, a positive constant $C_0$, and a sequence $\{\epsilon_n\}_{n\geq 1}$ of positive real numbers converging to 0
  such that for any $n\geq n_0$ and any $z\in \DD_{0.8 r_n}$, we have
\begin{equation*}%\label{equphiuniform}
d(\phi_n(z),\phi_{(x,\delta(s_{\infty}, v_b))}(z))\leq  d(s_n(b_n),x)+\epsilon_n|z|+C_0 \ell^{-1}_n|z|^2.
\end{equation*}
\end{lem}
\begin{proof}
We are going to consider suitable liftings of the sections via the covering map $\Lie(\CA_0)\to \CA_0$.
For every element $a\in \Lie(\CA_0)$, we still write $a$ for the constant map $\DD\to \Lie(\CA_0)$ of value $a$. 

Let $C^{\R\text{-}\an}(\D, \Lie(\CA_0))$ be the real vector space of real analytic maps $h:\DD\to \Lie(\CA_0)$, whose vector space structure is induced by that of $\Lie(\CA_0)$.  
Now we define an $\RR$-linear homomorphism 
$$\mathfrak b:V(A,K)\lra C^{\R\text{-}\an}(\D, \Lie(\CA_0)), \quad s\longmapsto s^{\flat}.$$
First, for every element $s\in A(K)$, let $s^*:\DD\to \Lie(\CA_0)$ be a (continuous) lifting of 
$\beta \circ s:\DD\to \CA_0$ and set $s^{\flat}:=s^*-s^{*}(0)$.  This definition does not depend on the choice of the lifting $s^*$.
Second, the map $s\mapsto s^{\flat}$ is a group homomorphism for $s\in A(K)$, whose kernel contains $A^{K/\CC}(\CC)$, so it induces the homomorphism 
$\mathfrak b$ on $V(A,K)$.

By Theorem \ref{non-deg2}, the homomorphism $\mathfrak b$ is injective.
By definition, for every $s\in V(A,K)$, the differential map 
$$(\d s^{\flat})_t: T_t\D\lra \Lie(\CA_0),$$
at $t\in \D$ is the composition
$$
T_t\D \stackrel{\delta_t(s,\cdot)}{\lra} \Lie(\CA_t) 
\lra T_{(t,0)}(\CA) \stackrel{\d \beta}{\lra} \Lie(\CA_0).
$$
Here the second arrow is induced by the immersion $\CA_t\to \CA$, and $(t,0)$ represents the image of $0$ under this immersion. 
%$$(\d s^{\flat})_t=(\d \beta)_{t,0}\circ \delta_t(s,\cdot): T_t\D\lra \Lie(\CA_0),$$
We view $(\d s^{\flat})_t$ as an $\R$-linear morphism from $\C=T_t\D$  to $\Lie(\CA_0),$ which is the linear part of $s^{\flat}$ at $t.$
Define the total differential map 
$$(s^{\flat})': \D\lra \Hom_{\R}(\C,\Lie(\CA_0)), \quad
t \longmapsto (\d s^{\flat})_t.$$
Similarly, define  
$$(s^{\flat})'': \D\lra \Hom_{\R}(\C,\Hom_{\R}(\C,\Lie(\CA_0))), \quad
t \longmapsto (d(s^{\flat})')_t.$$
Both $(s^{\flat})'$ and $(s^{\flat})''$ are real analytic and thus continuous.

Fix an euclidean norm $|\cdot|$ on $V(A,K)$.  
Still denote by $|\cdot|$ the norms on $\Hom_{\R}(\C,\Lie(\CA_0))$ and $\Hom_{\R}(\C,\Hom_{\R}(\C,\Lie(\CA_0)))$ induced by the standard norm on $\C$ and the norm $|\cdot|$ on  $\Lie(\CA_0).$
Since $V(A,K)$ is a finite-dimensional $\R$-space, there is a constant $C>0$ such that $$\max\{|(s^{\flat})''(t)|,|(s^{\flat})'(t)|\}\leq C|s|$$ for every $s\in V(A,K)$ and $|t|\leq 0.9.$
Then both $s^{\flat}$ and $(s^{\flat})'$ are $C|s|$-Lipschitz on $\overline{\DD(0,0.9)}$ for every $s\in V(A,K)$.
%\textcolor{red}{(Note: still problem of derivative.)}

Set $x_n:=\beta(s_n(b_n))\in \CA_0.$
Let $x^\dagger$ be a lifting of $x$ in $\Lie(\CA_0)$, and $x_n^\dagger$ be a lifting of $x_n$ in $\Lie(\CA_0)$ for every $n\geq 1$.
Since $x_n\to  x$, we may assume $x_n^{\dagger}\to x^{\dagger}.$ 
Set 
$$s_n^\dagger:=x_n^\dagger+s_n^\flat-s_n^{\flat}(b_n),
\qquad s_{\infty}^\dagger:=x^{\dagger}+s_{\infty}^\flat.$$
Then $s_n^\dagger$ is a lifting of $\beta\circ s_n$ by $s_n^\dagger(b_n)=x_n^\dagger$. 
By condition (2), $\ell_n^{-1}{s_n}$ converges to $s_\infty$. 
We may write 
$$s_n=\ell_ns_\infty+u_n, \quad |u_n|/\ell_n\to 0.$$
%Since $\d s_{\infty}^\flat(0)=\delta(s_{\infty}, v_b)$, we may write 
%$s_{\infty}^\flat(t)=\delta(s_{\infty}, v_b)t+w(t).$ 
%There is $C_0>0$ such that $|w(t)|\leq C_0|t|^2$ for $|t|\leq 1/2.$
Denote 
$$\phi_n^{\dagger}:=s_n^{\dagger}\circ p_n: \DD_{r_n}\lra \Lie(\CA_0).$$
 It is a lifting of $\phi_n.$ 
Then for every $z\in \D_{r_n}$, we have 
\begin{equation*}%\label{equasndagger}
\begin{split}
\phi_n^{\dagger}(z)
=&\, s_n^{\dagger}(b_n+\ell^{-1}_nz)\\
=&\, x_n^\dagger+s_n^\flat(b_n+\ell^{-1}_nz)-s_n^{\flat}(b_n) \\
=&\, x_n^{\dagger}+\ell_n(s_{\infty}^\flat(b_n+\ell^{-1}_nz)-s_{\infty}^\flat(b_n))+(u_n^{\flat}(b_n+\ell^{-1}_nz)-u_n^\flat(b_n)).
\end{split}
\end{equation*}
It follows that
\begin{equation}\label{equphin}
\begin{split}
|\phi_n^{\dagger}(z)-(x^{\dagger}+\delta(s_{\infty}, v_b)z)|\leq & |x_n^{\dagger}-x^{\dagger}|+|u_n^{\flat}(b_n+\ell^{-1}_nz)-u_n^\flat(b_n)| 
+\\ 
&|\ell_n(s_{\infty}^\flat(b_n+\ell^{-1}_nz)-s_{\infty}^\flat(b_n))-\delta(s_{\infty}, v_b)z|.
\end{split}
\end{equation}

Since $b_n\to 0$, there exists a positive integer $n_0$ such that $|b_n|<0.1$ for every $n\geq n_0$.
Then for $n\geq n_0$ and for $z\in \DD_{0.8r_n}$, 
$$
|b_n+\ell^{-1}_nz|< 0.1+ 0.8 r_n/\ell_n < 0.9.
$$ 
Thus the above Lipschitz condition holds for $b_n+\ell^{-1}_nz$ and $b_n$.
It follows that
\begin{equation}\label{equun}|u_n^{\flat}(b_n+\ell^{-1}_nz)-u_n^\flat(b_n)|\leq C|u_n|\cdot|\ell^{-1}_nz|.
\end{equation}

On the other hand, since 
$$s_{\infty}^\flat(b_n+\ell^{-1}_nz)-s_{\infty}^\flat(b_n)=((s_{\infty}^\flat)'(b_n))(\ell^{-1}_nz)+O(|\ell^{-1}_nz|^2),$$
 there is $C_0>0$ such that 
\begin{equation}\label{equsinfty}
|(s_{\infty}^\flat(b_n+\ell^{-1}_nz)-s_{\infty}^\flat(b_n))-((s_{\infty}^\flat)'(b_n))(\ell^{-1}_nz)|\leq C_0|\ell^{-1}_nz|^2.
\end{equation}
Since $(s^{\flat}_{\infty})'$ is $C|s_{\infty}|$-Lipschitz and $\delta(s_{\infty}, \cdot)=((s_{\infty}^\flat)'(0))(\cdot),$
\begin{equation}\label{equsinftyder}
|(s_{\infty}^\flat)'(b_n)-\delta(s_{\infty}, \cdot)|\leq C|s_{\infty}||b_n|.
\end{equation}
Note that on the $\C$-vector space $\Lie\CA_0$, we have $\delta(s_{\infty}, z)=\delta(s_{\infty}, v_b)z$.
%\textcolor{red}{(Note: The problem is still that $s_{\infty}^\flat$ is not holomorphic, but the induced map on the tangent space is $\CC$-linear.)}
Combining inequalities (\ref{equphin}), (\ref{equun}),(\ref{equsinfty}) and (\ref{equsinftyder}), we have 
\begin{equation} \label{equfinal}
|\phi_n^{\dagger}(z)-(x^{\dagger}+\delta(s_{\infty}, v_b)z)|\leq 
|x_n^{\dagger}-x^{\dagger}|+C|z|\cdot|u_n|/l_n+C_0\ell^{-1}_n|z|^2+C|s_{\infty}||b_nz|.
\end{equation}
By definition, 
$$
d(\phi_n(z),\phi_{(x,\delta(s_{\infty}, v_b))}(z))
\leq |\phi_n^{\dagger}(z)-(x^{\dagger}+\delta(s_{\infty}, v_b)z)|+|b_n|.
$$
Moreover, 
$$
d(s_n(b_n),x)=|x_n^{\dagger}-x^{\dagger}|+|b_n|,
$$
which holds for sufficiently large $n$ as $x_n^{\dagger}\to x^{\dagger}$. 
Then the lemma is a consequence of (\ref{equfinal}) by the two distance relations. 
\end{proof}

\subsubsection*{Proof of Theorem \ref{linear entire curve}(a)}
Since the on $\CA_{\DD_{0.9}}$ induced by the distance function $d$ is exactly the one induced by the K\"ahler metric $\alpha_{\CA}$,  Part (a) is a consequence of Lemma \ref{error}.

%
%Part (a) is a consequence of Lemma \ref{error} by the claim that the distance function $d$ on $\CA_{\DD_{0.9}}$ used in the lemma is dominated by a positive constant multiple of the distance function induced by the K\"ahler metric
%$\alpha_{\CA}$.
%
%To prove the claim, by compactness of $\CA_{\overline{\DD_{0.9}}}$, we can replace the K\"ahler metric $\alpha_{\CA}$
%by any Riemannian metric on $\CA_\DD$.
%In particular, by the Betti isomorphism $(\beta,\pi): \CA_\DD \to \CA_0\times \DD$, 
%we can replace $\alpha_{\CA}$ by a Riemannian metric of the form 
%$\alpha_{\CA}'=\beta^*\alpha_{\CA_0}+\pi^*\alpha_{\DD}$, where 
%$\alpha_{\CA_0}$ and $\alpha_{\DD}$ are Riemannian metrics on ${\CA_0}$ and ${\DD}$ respectively. 
%Take $\alpha_{\DD}$ to be the standard metric on $\DD$, and take 
%$\alpha_{\CA_0}$ to be the metric induced by the euclidean norm $|\cdot|$ on $\Lie(\CA_0)$.
%Then the distance function on $\CA_\DD\simeq \CA_0\times \DD$ with respect to $\alpha_{\CA}'$ is just 
%$$d'((x_1,t_1),(x_2,t_2))= \sqrt{d_{\CA_0}(x_1,x_2)^2+ |t_1-t_2|^2}.$$
%Compare it with 
%$$d((x_1,t_1),(x_2,t_2))= d_{\CA_0}(x_1,x_2)+ |t_1-t_2|.$$
%We simply have 
%$$
%d((x_1,t_1),(x_2,t_2))\leq \sqrt2 d'((x_1,t_1),(x_2,t_2)).
%$$
%This finishes the proof.

\subsubsection*{Proof of Theorem \ref{linear entire curve}(c)}

Let $\{r_n'\}_{n\geq1}$ be a sequence in the setting of Theorem \ref{linear entire curve}(a). 
Recall that Lemma \ref{error} asserts that for any $n\geq n_0$ and any $z\in \DD_{0.8 r_n'}$, 
\begin{equation*}%\label{equphiuniform}
d(\phi_n(z),\phi_{(x,\delta(s_{\infty}, v_b))}(z))\leq  d(s_n(b_n),x)+\epsilon_n|z|+C_0 \ell^{-1}_n|z|^2.
\end{equation*}
Let $\{r_n\}_{n\geq1}$ be a sequence of positive integers such that 
$$
r_n<0.8r_n',\quad 
r_n\lra \infty, \quad
\epsilon_n r_n\lra 0,\quad
\ell^{-1}_n r_n^2\lra 0. 
$$
Since $\max\{\epsilon_n, \ell^{-1/2}_n\}\to 0$, such $\{r_n\}_{n\geq 1}$ exists.
By this choice, for any $n\geq n_0$ and for any $z\in \DD_{r_n}$, 
\begin{equation}\label{equphiuniform}
d(\phi_n(z),\phi_{(x,\delta(s_{\infty}, v_b))}(z))\lra 0.
\end{equation}
The convergence is uniform in $z$. 
By this key property, we will prove that part (c) is satisfied by the sequence $\{r_n\}_{n\geq1}$.

For the sake of part (c), let $\psi:\C\to \CA$ be a limit of a re-parametrization $\{\psi_n\}_{n\geq 1}$ of $\{\phi_n\}_{n\geq 1}$. One may write
$\psi_n=\phi_n\circ q_n: \DD_{R_n}\to \sA,$
where 
$\{R_n\}_{n\geq 1}$ is a sequence of positive real numbers converging to infinity and $\{q_n: \DD_{R_n} \to \D_{r_n}\}_{n\geq 1}$ is a sequence of holomorphic maps.
After taking a subsequence, we may assume that $\psi_n$ converges to $\psi$ uniformly on every disk $\D_R$ for $R>0.$
By (\ref{equphiuniform}), for every $z\in \D_R$, we have 
\begin{equation*}%\label{equphiuniformpsi}
d(\phi_n(q_n(z)),\phi_{(x,\delta(s_{\infty}, v_b))}(q_n(z)))\lra 0.
\end{equation*}
As a consequence, $\phi_{(x,\delta(s_{\infty}, v_b))}\circ q_n$ converges to $\psi$ uniformly.

Set 
$$y:=\psi(0)\in \CA_0.$$ 
Then $\phi_{(x,\delta(s_{\infty}, v_b))}(q_n(0))\to y.$ It follows that $y\in \overline{\phi_{(x,\delta(s_{\infty}, v_b))}(\C)}.$
This gives the point $y$ in the theorem. 

For further properties, set 
$$h_n(z):=q_n(z)-q_n(0).$$ 
Then $h_n(0)=0$. For sufficiently large $n$ and for $z\in \D_{R_n}$, we have 
\begin{equation*}
\begin{split}
\phi_{(y,\delta(s_{\infty}, v_b))}\circ h_n&=y+
\phi_{(0,\delta(s_{\infty}, v_b))}(h_n)\\
& =(y-\phi_{(x,\delta(s_{\infty}, v_b))}(q_n(0)))+\phi_{(x,\delta(s_{\infty}, v_b))}\circ q_n.
\end{split}
\end{equation*}
Since $(y-\phi_{(x,\delta(s_{\infty}, v_b))}(q_n(0)))\to 0$, we have 
$\phi_{(y,\delta(s_{\infty}, v_b))}\circ h_n\to \psi$ uniformly.

We claim that $h_n$ uniformly converges to a morphism $h: \C\to \C$. In particular, we have $\psi=\phi_{(y,\delta(s_{\infty}, v_b))}\circ h$.
It implies that $\psi(\CC)\subseteq \phi_{(y,\delta(s_\infty,v_b))}(\CC)$. Since $\psi$ is non-constant, $h$ is non-constant. Since $\C\setminus h(\C)$ has at most 1 element, 
$\phi_{(y, \delta(s_\infty,v_b))}(\CC)\setminus \psi(\CC)$ has at most 1 element. We then get property (c).

It remains to prove the claim. Denote by $\wt y$ a lifting of $y$ in $\Lie(\CA_0)$.
Let $\wt\phi:\CC \to \Lie(A)$ be the unique lifting of 
$\phi_{(y,\delta(s_\infty,v_b))}$ satisfying $\wt\phi(0)=\wt y$.
There is $c>0$ such that the restriction of the map $\Lie(\CA_0)\to \CA_0$ is injective on the closed ball $\{u\in \Lie(\CA_0)|\,\, |u|\leq c\}.$
So for every $u\in \Lie(\CA_0)$ and $v\in \CA_0$, there is at most one lifting $\wt v$ of $v$ satisfying $|u-\wt v|\leq c$.

Since $\phi_{(y,\delta(s_{\infty}, v_b))}\circ h_n\to \phi_{(y,\delta(s_{\infty}, v_b))}\circ h$ uniformly, for every $\epsilon\in (0,c/2)$ there is $N\geq 1$, such that for every $n\geq N$ and $z\in \D_{r_n},$ $$d_{\CA_0}(\phi_{(y,\delta(s_{\infty}, v_b))}\circ h_n, \phi_{(y,\delta(s_{\infty}, v_b))}\circ h)<\epsilon.$$
So for every $z\in \D_{r_n}$, 
$$|\wt\phi\circ h(z)-\wt\phi\circ h_n(z)|\in [0,\epsilon)\cup (c,+\infty).$$
Since $\D_{r_n}$ is connected and $|\wt\phi\circ h(0)-\wt\phi\circ h_n(0)|=0$, we have 
$|\wt\phi\circ h(z)-\wt\phi\circ h_n(z)|< \epsilon$ for every $z\in \D_{r_n}.$
Hence $|h(z)-h_n(z)|< \epsilon/|\delta(s_{\infty}, v_b)|$ for every $z\in \D_{r_n},$ which proves the claim.

\subsection{Locations of limit points}  \label{sec location}

Another key result to prove Theorem \ref{main} is the following theorem about  locations of limit points, which is closely related to locations of entire curves by combining Theorem \ref{linear entire curve}. 

\begin{thm} \label{limit point}
Let $B$ be a smooth quasi-projective curve over $\CC$, and $\pi:\CA\to B$ be an abelian scheme.
Let $K=\CC(B)$ be the function field of $B$, and $A$ be the generic fiber of $\pi:\CA\to B$. 
Assume that the $K/\CC$-trace of $A$ is 0, and that the natural map $\End_{K}(A)\to \End_{\ol K}(A_{\ol K})$ is an isomorphism.
Let $Z\subsetneq A$ be a Zariski closed subset of $A$, and $\CZ$ be the Zariski closure of $Z$ in $\CA$.   
Let $\{s_n\}_{n\geq1}$ be an infinite sequence of distinct elements in $(A\setminus Z)(K)$. 
Then there is a subset $S$ of $B(\CC)$ of measure 0 such that for every $b\in B(\CC)\setminus S$, there exists a point $y\in \CA_b\setminus \CZ_b$ which is a limit point of $\{s_n(B)\}_{n\geq1}$, i.e.
there is a sequence $\{b_n\}_{n\geq1}$ in $B(\CC)$ such that
the sequence $\{s_n(b_n)\}_{n\geq1}$ has a subsequence converging to $y$ in $\CA$.
\end{thm}

\begin{proof}

By the Lang--N\'eron theorem (cf. \cite[Thm. 3.4]{XY2023a}), the abelian group $A(K)$ is finitely generated. 
By replacing $\{s_n\}_{n\geq1}$ by a subsequence, 
we can assume that there is a sequence $\ell_n$ of positive real numbers converging to infinity such that $\ell_n^{-1}s_n$ 
converges to a nonzero element $s_\infty$ in $A(K)_\RR$. 
In other words, Theorem \ref{linear entire curve}(2) is satisfied.
By Theorem \ref{non-deg2}, up to removing a subset of $B(\CC)$ of measure 0, 
we can assume that Theorem \ref{linear entire curve}(3) is satisfied for $b\in B(\CC)$. 

Our guiding principle is that once Theorem \ref{linear entire curve}(2)(3) is satisfied for
 fixed $\{s_n\}_{n\geq1}$ and fixed $b\in B(\CC)$ (as assumed by us), the following are equivalent:
\begin{enumerate}[(i)]
\item There exists $y\in \CA_b\setminus \CZ_b$ which is a limit point of $\{s_n(B)\}_{n\geq1}$.
\item There exists $x\in \CA_b$ which is a limit point of $\{s_n(B)\}_{n\geq1}$ such that the linear leaf $\CF(x,\delta(s_\infty,v_b))$ is not contained in $\CZ_b$. Here $v_b\in T_bB$ is any nonzero vector.
\end{enumerate}
Note that (i) implies (ii) by setting $x=y$. 
Conversely, (ii) implies (i) by taking a point 
$y\in \CF(x,\delta(s_\infty,v_b))\setminus \CZ_b$.
In fact, by Theorem \ref{linear entire curve}(a), 
$\phi_{(x,\delta(s_\infty,v_b))}:\CC\to \CA$ is a limit of a re-parametrization of 
$\{s_n\}_{n\geq1}$, so the definition implies that every point of $ \CF(x,\delta(s_\infty,v_b))=\phi_{(x,\delta(s_\infty,v_b))}(\CC)$ is a limit point of $\{s_n(B)\}_{n\geq1}$.

Denote by $G$ the smallest abelian subvariety of $A$ such that 
$s_\infty$ lies in $G(K)_\RR$. 
Denote by $\CG$ the Zariski closure of $G$ in $\CA$. 
By Proposition  \ref{equality of zariski closure}, up to removing a subset of $B(\CC)$ of measure 0, 
we can assume that the Zariski closure of  
$\CF(0,\delta(s_{\infty},v_b))$ in $\CA_b$ is equal to $\CG_b$. 

If $G=A$, the result is easy. 
In fact, take $x\in\CA_b$ to be any limit point of $\{s_n(b)\}_{n\geq1}$, which exists by compactness.
Then the Zariski closure of  
$\CF(x,\delta(s_{\infty},v_b))$ in $\CA_b$ is equal to $\CG_b+x=\CA_b$. 
It follows that
$\CF(x,\delta(s_\infty,v_b))$ is not contained in $\CZ_b$. 
This proves the result by the equivalence of (i) and (ii). 

Now we prove the theorem for the general case by induction on the number of simple isogeny factors of $A$. 
If $A$ is simple, this is included in the case $G=A$.
Assume that $A$ is not simple and that $G\neq A$.

For $b\in B(\CC)$ as above, denote by $E_b$ the union of all subvarieties of $\CZ_b$ which are of the form $\CG_b+a$ for some $a\in \CA_b(\CC)$. 
If $E_b$ is empty, the result holds by the equivalence of (i) and (ii). 
In the following, we assume that $E_b$ is non-empty. 

We can realize $E_b(\CC)$ as the closed points of a subscheme $\CE_b$ of $\CZ_b$, and lift it to a subscheme $\CE$ of $\CZ$ over $B$. 
In fact, consider the composition 
$$
\CZ \lra \CA \lra \CA/\CG.
$$
By semicontinuity of dimensions of fibers, there is a maximal closed subset $\CE$ of $\CZ$ such that every non-empty fiber of $\CE\to \CA/\CG$ has a dimension equal to $\dim G$.
Then the closed points $\CE_b(\CC)$ of the fiber $\CE_b$ is exactly $E_b(\CC)$.  

Now we take the quotient $(\CA', \CZ')=(\CA/\CG, \CE/\CG)$ over $B$.
Denote by $(A',Z')=(\CA'_K,\CZ'_K)$ the generic fibers. 
By removing finitely many $b\in B(\CC)$, we can assume that $\CE$ is flat at $b\in B(\CC)$, and that $\CE_b\to \CZ'_b$ is a quotient by $\CG_b$. 

We plan to apply the induction hypothesis to $(A', Z')$.
Denote by $s_n'\in \CA'(B)$ the image of $s_n\in \CA(B)$. 
The following discussion is based on whether $\{s_n'\}_{n\geq1}$ has infinitely many distinct terms.

If $\{s_n'\}_{n\geq1}$ has infinitely many distinct terms,  by induction, there is a closed point $x'\in \CA'_b\setminus \CZ'_b$ which is a limit point of $\{s_n'(b_n)\}_{n\geq1}$ for a sequence $\{b_n\}_{n\geq1}$ of $B$. 
By compactness, $\{s_n(b_n)\}_{n\geq1}$ has a limit point $x\in \CA_b$ whose image in $\CA'_b$ is $x'$.
Recall that the Zariski closure of $\CF(x,\delta(s_\infty,v_b))$ in $\CA_b$ is just $\CG_b+x$.
By $x'\in \CA'_b\setminus \CZ'_b$, 
we have $x\in \CA_b\setminus \CE_b$.  
Then $\CG_b+x$ is not contained in $\CE_b$, and thus not contained in $\CZ_b$.
Then $\CF(x,\delta(s_\infty,v_b))$ is not contained in $\CZ_b$. 
This proves the result by the equivalence of (i) and (ii). 

If $\{s_n'\}_{n\geq1}$ has only finitely many distinct terms, then by passing to a subsequence, we can assume that all $s_n'$ are equal in $\CA'(B)$.
Then $\{s_{n}\}_{n\geq1}$, viewed as a sequence in $A(K)$, lies in a fiber of $A\to A'$.  
This fiber is equal to the translate $G+s_{1}$ of $G$ in $A$. 
As $s_1$ is not contained in $Z$ by assumption, we see that $G+s_1$ is not contained in $Z$.
By removing finitely many $b\in B(\CC)$, we can assume that 
 the specialization $\CG_b+s_1(b)$ is not contained in $\CZ_b$.
Then $\CG_b+s_1(b)$ does not intersect $\CE_b$, since $\CE_b$ is a union of cosets of $\CG_b$.  
Take an arbitrary limit point $x$ of $\{s_n(b)\}_{n\geq 1}$ in $\CG_b+s_1(b)$.
Then $\CF(x,\delta(s_\infty,v_b))$ is not contained in $\CZ_b$, as its Zariski closure is $\CG_b+s_1(b)$.
The result follows by the equivalence of (i) and (ii). 
This finishes the proof.
\end{proof}

\subsection{Proof of the main theorem}  \label{sec proof}

In this subsection, we prove Theorem \ref{main}. The theorem is restated as follows.

\begin{thm} [Theorem \ref{main}] \label{main copy}
Let $K$ be a finitely generated field over a field $k$ of characteristic 0. 
Let $X$ be a projective variety over $K$ with a finite morphism $f:X\to A$ for an abelian variety $A$ over $K$. 
Assume that the $\ol K/\bar k$-trace of $A_{\ol K}$ is 0. 
Then $(X\setminus \Sp_{\rm alg}(X))(K')$ is finite for any finite extension $K'$ of $K$.
\end{thm}

\subsubsection*{Base fields} 

Recall that Theorem \ref{main copy} is stated for a finite generated extension $K/k$ of characteristic 0.
As in \cite[\S4.2]{XY2023a}, we can reduce it to the case that $\trdeg(K/k)=1$ and $k=\CC$. 

We first reduce the theorem from $\trdeg(K/k)>1$ to $\trdeg(K/k)=1$.
Assume $\trdeg(K/k)>1$.
Let $k_1$ be an intermediate field of $K/k$ such that 
$\trdeg(K/k)=1$.
We hope that the result for $(K/k,X,A)$ is implied by that for  
$(K/k_1,X,A)$.
For that, it suffices to find such $k_1$ such that 
the $\ol K/\bar k_1$-trace of $A$ is 0, under the condition that
the $\ol K/\bar k$-trace of $A$ is 0. 
The existence of $k_1$ is a consequence of \cite[Cor. 3.5(1)]{XY22}. 
In fact, denote by $A_1,\dots, A_r$ the simple isogeny components of $A$ over $K$.
We only need to find $k_1$ such that every $A_i$ is not defined over $k_1$, but this means that $k_1$ is not contained in any $k_{A_i}$ obtained by the loc. cit.. 

Now we reduce the theorem to $k=\CC$ (assuming $\trdeg(K/k)=1$).
We can assume that $k$ is algebraically closed in $K$.
For the sake of contradiction, assume that $(X\setminus \Sp_{\rm alg}(X))(K)$ is infinite. Then it contains a countable subset $\Sigma$.
As in \cite[\S4.2]{XY2023a}, by the Lefschetz principle, we can descend 
$(K, f:X\to A, \Sigma)$ to a finitely generated field $k_0$ over $\QQ$, and then take the base change to $\CC$ via an embedding $k_0\to \CC$. 
In this process, the truth of the theorem does not change due to two properties. 
First, taking special set of projective varieties is 
stable under base change by
\cite[Lem. 4.2]{XY2023a}.
Second, taking the $\bar K/\bar k$-trace of $A$ is stable under base change by \cite[Thm. 6.8]{Conrad}. 

By replacing $K$ by a finite extension if necessary, we can further assume that the natural map $\End_{K}(A)\to \End_{\ol K}(A_{\ol K})$ is an isomorphism. We can further assume that $K'=K$. 
Then we have reduced the theorem to the following special case.

\begin{thm}  \label{main variant}
Let $K$ be a function field of one variable over $\CC$.
Let $X$ be a projective variety with a finite morphism $f:X\to A$ for an abelian variety $A$ over $K$. 
Assume that the $K/\CC$-trace of $A$ is 0, and that the natural map $\End_{K}(A)\to \End_{\ol K}(A_{\ol K})$ is an isomorphism.
Then $(X\setminus \Sp_{\rm alg}(X))(K)$ is finite.
\end{thm}

\subsubsection*{Subvarieties of abelian varieties} 

We first prove Theorem \ref{main variant} 
assuming that $f:X\to A$ is a closed immersion. 
This is an easy consequence of Theorem \ref{limit point} and 
Theorem \ref{linear entire curve}.

In fact, let $B$ be a smooth curve over $\CC$ with function field $K$, and let $\pi:\CA\to B$ be an abelian scheme extending the abelian variety $A\to \Spec K$. 
Let $\CX$ be the Zariski closure of $X$ in $\CA$.
Assume that $(X\setminus \Sp_{\rm alg}(X))(K)$ is infinite, so it has an infinite sequence $\{s_n\}_{n\geq1}$ of distinct terms.
Apply Theorem \ref{limit point} to $Z=\Sp_{\rm alg}(X)$ and $\{s_n\}_{n\geq1}$. 
We obtain a point $b\in B(\CC)$ together with a limit point 
$x\in \CA_b\setminus \Sp_{\rm alg}(\CX_b)$ of 
$\{s_n(B)\}_{n\geq1}$.
Here we assume that $\Sp_{\rm alg}(\CX_b)$ is the specialization of the Zariski closure of $\Sp_{\rm alg}(X)$ by \cite[Lem. 4.4, Lem. 4.5]{XY2023a}.
Apply Theorem \ref{linear entire curve} to $\{s_n\}_{n\geq1}$ and the limit point $x\in \CA_b\setminus \Sp_{\rm alg}(\CX_b)$.
We obtain a linear entire curve 
$\phi_{(x,\delta(s_\infty,v_b))}:\CC\to \CA_b$.
As $s_n(B)$ lies in $\CX$, the entire curve actually lies in $\CX_b$.
As $x\notin \Sp_{\rm alg}(\CX_b)$, we see that $\phi_{(x,\delta(s_\infty,v_b))}(\CC)$ is not contained in $\Sp_{\rm alg}(\CX_b)$.
By \cite[Thm. 4.1]{XY2023a}, which is essentially due to Kawamata \cite{Kawamata81} and Yamanoi \cite{Yamanoi15}, $\Sp_{\rm alg}(\CX_b)=\Sp_{\rm an}(\CX_b)$. 
It contradicts to the definition of $\Sp_{\rm an}(\CX_b)$. 
This finishes the proof assuming that $f:X\to A$ is a closed immersion.

\subsubsection*{Weaker statement} 

In Theorem \ref{main variant}, consider the weaker statement that $(X\setminus f^{-1}(f(\Sp_{\rm alg}(X)))(K)$ is finite.
Here we prove that the weaker statement for all pairs $(X,f)$
(with fixed $(A,K)$)  as in the theorem implies the original stronger statement. 

Induct on $\dim \Sp_{\rm alg}(X)$.
The results are equivalent if $\Sp_{\rm alg}(X)$ is empty, and note that $\dim \Sp_{\rm alg}(X)\neq 0$ by definition. 
If $\dim \Sp_{\rm alg}(X)\geq 1$, let $Y_1,\dots, Y_r$ be the irreducible components of $f^{-1}(f(\Sp_{\rm alg}(X))$ which are not contained in $\Sp_{\rm alg}(X)$. 
Note that $\Sp_{\rm alg}(Y_i)\subseteq \Sp_{\rm alg}(X)\cap Y_i$ has a dimension strictly smaller than $\Sp_{\rm alg}(X)$. 
By induction, $(Y_i\setminus \Sp_{\rm alg}(Y_i))(K)$ is finite, so $(Y_i\setminus \Sp_{\rm alg}(X))(K)$ is finite. 
It follows that
$$(X\setminus \Sp_{\rm alg}(X))(K)=(X\setminus f^{-1}(f(\Sp_{\rm alg}(X)))(K)\
\cup\ \big(\cup_{i=1}^r(Y_i\setminus \Sp_{\rm alg}(X))(K)\big)$$
is finite.

\subsubsection*{Proof of Theorem \ref{main variant}: entire curves} 

Now we prove Theorem \ref{main variant}. 
The proof is based on Theorem \ref{limit point} and 
Theorem \ref{linear entire curve}, and also uses ideas of the proof of Theorem \ref{limit point}.

As above, it suffices to prove that $(X\setminus f^{-1}(f(\Sp_{\rm alg}(X)))(K)$ is finite. 
Assume that this does not hold.
So there is an infinite sequence $\{x_n\}_{n\geq 1}$ of distinct elements of $(X\setminus f^{-1}(f(\Sp_{\rm alg}(X)))(K)$.
Denote $s_n=f(x_n)$. 
Then we have an infinite sequence $\{s_n\}_{n\geq 1}$ in $(A\setminus f(\Sp_{\rm alg}(X))(K)$.
By passing to a subsequence, we can assume that  
$\{s_n\}_{n\geq 1}$ has distinct terms.

Let $B$ be a smooth curve over $\CC$ with function field $K$.
Let $\pi:\CA\to B$ be an abelian scheme extending the abelian variety $A\to \Spec K$, and let $\tau:\CX\to B$ be a (projective and flat) integral model of $X$ over $B$. 
Assume that $f:X\to A$ extends to a finite morphism $f:\CX\to \CA$. 
All these exist by replacing $B$ by an open subscheme if necessary.
Denote $Z=f(\Sp_{\rm alg}(X))$, which is a Zariski closed subset of $A$.
Denote by $\overline{\Sp_{\rm alg}(X)}^{\zar}$ the Zariski closure of $\Sp_{\rm alg}(X)$ in $\CX$, and thus $\CZ=f(\overline{\Sp_{\rm alg}(X)}^{\zar})$ is the Zariski closure of $Z$ in $\CA$.

Up to replacing $\{s_n\}_{n\geq1}$ by a subsequence, 
we can assume that there is a sequence 
$\{\ell_{n}\}_{n\geq1}$ of positive real numbers converging to infinity such that the image of $\ell_{n}^{-1} s_{n}$ in 
$A(K)_\RR$ converges to a nonzero element $s_{\infty}\in A(K)_\RR$.
Apply Theorem \ref{limit point} to $Z=f(\Sp_{\rm alg}(X))$ and $\{s_n\}_{n\geq1}$. 
We obtain a point $b\in B(\CC)$ together with a limit point 
$x\in \CA_b\setminus \CZ_b$ of 
$\{s_n(b_n)\}_{n\geq1}$ for a sequence $\{b_n\}_{n\geq1}$ of $B(\CC)$.
By Theorem \ref{non-deg2}, 
we can further assume that the transfer map $\delta(s_\infty,\cdot):T_{b}B\to \Lie(\CA_{b})$ is nonzero. 
Apply Theorem \ref{linear entire curve} to $\{s_n\}_{n\geq1}$ and the limit point $x\in \CA_b\setminus \CZ_b$ of $\{s_n(b_n)\}_{n\geq1}$.
We obtain a linear entire curve 
$\phi_{(x,\delta(s_\infty,v_b))}:\CC\to \CA_b$,
which is the unique limit of the re-parametrization $\{\phi_n:\DD_{r_n}\to \CA\}_{n\geq1}$ of $\{s_n:B\to \CA\}_n$.

Now we want to ``lift'' the entire curve 
$\phi_{(x,\delta(s_\infty,v_b))}:\CC\to \CA$
to an entire curve on $\CX$.
Recall from Theorem \ref{linear entire curve} that $\{\phi_n:\DD_{r_n}\to \CA\}_{n\geq1}$ is defined by 
$$
\phi_n: \DD_{r_n}\lra \CA, \quad z\longmapsto s_n(b_n+\ell_n^{-1}z)
.$$ 
We may assume further that $\{r_n\}_{n\geq 1}$ satisfies the condition in Theorem \ref{linear entire curve}(c).
Then we define a sequence
$\{\tilde\phi_n:\DD_{r_n}\to \CX\}_{n\geq1}$ by 
$$
\tilde\phi_n: \DD_{r_n}\lra \CX, \quad z\longmapsto x_n(b_n+\ell_n^{-1}z)
$$ 
By definition, $\phi_n=f\circ \tilde\phi_n$. 
We plan to apply the Brody lemma (cf. \cite[Thm. 1.6]{XY2023a})
to the sequence $\{\tilde\phi_n:\DD_{r_n}\to \CX\}_{n\geq1}$
to produce an entire curve on $\CX_b$. 
Note that the domain of $\tilde\phi_n$ is $\DD_{r_n}$, so we modify it to $\DD$ by multiplying by $r_n$.
More precisely, we define a sequence
$\{\tilde\phi_n^\circ:\DD\to \CX\}_{n\geq1}$ by 
$\tilde\phi_n^\circ(z)=\tilde\phi_n(r_n z).$

We first check that $\{\tilde\phi_n^\circ:\DD\to \CX\}_{n\geq1}$ satisfies the growth condition in the Brody lemma.
Namely, we claim that 
$\|(\d \tilde\phi_n^\circ)_{0}(v_{\rm st})\|_{\alpha_{\CX}}$
converges to infinity, 
where $\alpha_{\CX}$ is a K\"ahler metric on a compactification of $\CX$.
In fact, by Theorem \ref{linear entire curve}(b), 
$\inf_n\|(\d \phi_n)_{0}(v_{\rm st})\|_{\alpha_{\CA}}>0$
for some K\"ahler metric $\alpha_{\CA}$ on $\CA$. 
By \cite[Lem. 2.1]{XY2023a}, there is a real constant $c>0$ such that $\alpha_{\CX}|_{\CX_\DD}-c f^*\alpha_{\CA}|_{\CX_\DD}$ is positive on $\CX_\DD$.
Then 
$\inf_n\|(\d \tilde\phi_n)_{0}(v_{\rm st})\|_{\alpha_{\CX}}>0$.
By $(\d \tilde\phi_n^\circ)_{0}=r_n(\d \tilde\phi_n)_{0}$, the growth condition is satisfied.

Finally, we can apply the Brody lemma (cf. \cite[Thm. 1.6]{XY2023a})
to the sequence $\{\tilde\phi_n^\circ:\D\to \CX\}_{n\geq1}$. 
Then some subsequence of some re-parametrization of $\{\tilde\phi_n^{\circ}:\DD\to \CX\}_{n\geq1}$
converges to an entire curve $\tilde\phi:\CC\to \CX$.
By the map $f:X\to A$, we have an entire curve  
$\psi=f\circ \tilde\phi:\CC\to \CA$.
Note that any re-parametrization of $\{\tilde\phi_n^{\circ}:\DD\to \CX\}_{n\geq1}$ is also a re-parametrization of $\{\tilde\phi_n:\DD_{r_n}\to \CX\}_{n\geq1}$, and $f:\CX\to \CA$ maps it to a re-parametrization of $\{\phi_n:\DD_{r_n}\to \CA\}_{n\geq1}$.
It follows that $\psi:\CC\to \CA$ is a limit of a re-parametrization of 
$\{\phi_n:\DD_{r_n}\to \CX\}_{n\geq1}$.

It is not reasonable to expect $\psi(\CC)=\phi_{(x, \delta(s_\infty,v_b))}(\CC)$ in general.
However, by Theorem \ref{linear entire curve}(c), $x\in \overline{\CF(x,\delta(s_\infty,v_b)))}=\overline{\psi(\CC)}.$
 Since $x\in \CA_b\setminus \CZ_b$,  $\psi(\CC)\not\subseteq \CZ_b$.
Thus
$\tilde\phi(\CC)\not\subseteq f^{-1}(\CZ_b)$.
By definition, $\CZ=f(\overline{\Sp_{\rm alg}(X)}^{\zar})$. 
We can assume that $\CZ_b=f(\overline{\Sp_{\rm alg}(X)}^{\zar}_b)$ by removing finitely many $b\in B(\CC)$. 
By \cite[Lem. 4.4, Lem. 4.5]{XY2023a}, we can assume that $\overline{\Sp_{\rm alg}(X)}^{\zar}_b= \Sp_{\rm alg}(\CX_b)$ by choosing $b\in B(\CC)$ suitably. 
It follows that
$f^{-1}(\CZ_b)=f^{-1}(f(\Sp_{\rm alg}(\CX_b))) \supseteq \Sp_{\rm alg}(\CX_b)$. 
It follows that 
$\tilde\phi(\CC)\not\subseteq \Sp_{\rm alg}(\CX_b)$.
By \cite[Thm. 4.1]{XY2023a}, which is essentially due to Kawamata \cite{Kawamata81} and Yamanoi \cite{Yamanoi15},
$\Sp_{\rm alg}(\CX_b)=\Sp_{\rm an}(\CX_b)$. 
This is a contradiction.

%More substantially, by Theorem \ref{linear entire curve}(d), the image
%$\psi(\CC)$ is actually contained in $\CG_b+x$, where $\CG$ is the smallest abelian subscheme of $\CA$ such that $\CG(B)_\RR$ contains $s_\infty$.
%By Lemma \ref{equality of zariski closure}, we can assume that 
%$\CG_b$ is the Zariski closure of $\phi_{(0, \delta(s_\infty,v_b))}(\CC)$. 
%It follows that the Zariski closures of 
%$\psi(\CC)$ and $\phi_{(x, \delta(s_\infty,v_b))}(\CC)$ are both equal to 
%$\CG_b+x$. 
%
%
%Hence, $\psi(\CC)\not\subseteq \CZ_b$ and thus
%$\tilde\phi(\CC)\not\subseteq f^{-1}(\CZ_b)$.
%By definition, $\CZ=f(\overline{\Sp_{\rm alg}(X)})$. 
%We can assume that $\CZ_b=f(\overline{\Sp_{\rm alg}(X)}_b)$ by removing finitely many $b\in B(\CC)$. 
%By \cite[Lem. 4.4, Lem. 4.5]{XY2023a}, we can assume that $\overline{\Sp_{\rm alg}(X)}_b= \Sp_{\rm alg}(\CX_b)$ by choosing $b\in B(\CC)$ suitably. 
%It follows that
%$f^{-1}(\CZ_b)=f^{-1}(f(\Sp_{\rm alg}(\CX_b))) \supseteq \Sp_{\rm alg}(\CX_b)$. 
%It follows that 
%$\tilde\phi(\CC)\not\subseteq \Sp_{\rm alg}(\CX_b)$.
%By \cite[Thm. 4.1]{XY2023a}, which is essentially due to Kawamata \cite{Kawamata81} and Yamanoi \cite{Yamanoi15},
%$\Sp_{\rm alg}(\CX_b)=\Sp_{\rm an}(\CX_b)$. 
%This is a contradiction.

%%%%%%%%%
\bibliographystyle{alpha}
\bibliography{dd}

%%%%%%%%%
\medskip
\noindent \small{Address: \textit{Beijing International Center for Mathematical Research, Peking University, Beijing 100871, China}}

\noindent \small{Email: \textit{xiejunyi@bicmr.pku.edu.cn}}

\medskip
\noindent \small{Address: \textit{Beijing International Center for Mathematical Research, Peking University, Beijing 100871, China}}

\noindent \small{Email: \textit{yxy@bicmr.pku.edu.cn}}

%\newpage
%
%(1). I still want to change the title, since we treat little about varieties of maximal albanese dimension, and the statement for that case is also not neat enough. 
%Some other options are:
%
%The geometric Bombieri--Lang conjecture for ramified covers of abelian varieties.
%
%The geometric Bombieri--Lang conjecture for varieties finite over abelian varieties.
%
%Partial heights and the geometric Bombieri--Lang conjecture.
%
%I like the third name, which emphasizes our new notion of partial heights. 
%
%\
%
%(2). What can we say in Theorem \ref{partial to BL0} when $X$ is singular? 
%The strong statement for the normalization of $X$ does not seem correct, so the result is weaker than that in Theorem \ref{hyperbolic main}. 
%Need to think about the example given by Xu.
%Need to find the correct statement...
%
%\
%
%(3). change $\Sp_{\rm alg}(X), \Sp_{\rm an}(X)$ to $\Sp_{\rm alg}(X), \Sp_{\rm an}(X)$?
%
%\
%
%(4). algebraically hyperbolic vs. algebraically hyperbolic?
%
%\

%------------------------------------------------------
\end{document}